\documentclass[12pt,a4paper]{article}
\usepackage{amsfonts,amssymb,amsmath,amsbsy}
\usepackage{latexsym}
\usepackage{graphicx}

\newtheorem{theo}{Theorem}
\newtheorem{lemma}[theo]{Lemma}
\newtheorem{prop}[theo]{Proposition}

\newtheorem{conj}[theo]{Conjecture}

\newcommand{\proof}{\noindent{\it Proof: }}
\newcommand{\proofbox}{\hfill \mbox{ $\Box$}\\}

\newcommand{\R}{{\mathbb R}}

\newcommand{\Z}{{\mathbb Z}}
\newcommand{\N}{{\mathbb N}}

\title{A note on triangulations of sum sets}
\author{K.J. B\"or\"oczky \and B. Hoffman}

\begin{document}
\maketitle
\begin{abstract}
In $\R^2$, for finite sets $A$ and $B$, we write $A+B=\{a+b:a\in A, b\in B\}$. We write ${\rm tr}(A)$ to denote the common number of triangles in any triangulation of the convex hull of $A$ using the points of $A$ as vertices. We consider the conjecture that ${\rm tr}(A+B)^{\frac12}\geq {\rm tr}(A)^{\frac12}+{\rm tr}(B)^{\frac12}$. If true, this conjecture would be a discrete, two-dimensional analogue to the Brunn-Minkowski inequality. We prove the conjecture in three special cases.
\end{abstract}
\section{Introdution}

In this paper, we write  $A,B$ to denote finite
subsets of $\R^d$,
and $|\cdot|$ to stand for their cardinality.
For objects $X_1,\ldots,X_k$ in $\R^2$,
$[X_1,\ldots,X_k]$ denotes their convex hull.
Our starting point is two classical results.
One is due to Freiman from the 1960's;
namely,
\begin{equation}
\label{linearsum}
|A+ B|\geq |A|+|B|-1,
\end{equation}
with equality if and only
if $A$ and $B$ are arithmetic progressions of the same difference.
The other result,
the Brunn-Minkowski inequality dates back to the 19th century.
It says that if $X,Y\subset \R^d$ are compact
sets then
$$
\lambda(X+Y)^{\frac1d}\geq \lambda(X)^{\frac1d}+\lambda(Y)^{\frac1d}
$$
where $\lambda$ stand for the Lebesgue measure,
and equality holds if $X$ and $Y$ are convex
homothetic sets.
This theorem has been successfully applied to
estimating the size of a a sumset say by Ruzsa, Green, Tao.
In turn various discrete analogues
of the Brunn-Minkowski inequality have been
established in papers by
Bollob\'as-Leader, Gardner-Gronchi, Green-Tao
and most recently by Grynkiewicz-Serra in the planar case.
All these papers use the method
of compression, which changes a
finite set into a set better suited for
sumset estimates, but cannot control the convex hull. See
G.A. Freiman \cite{Fre73} and  \cite{Fre02} for the earlier history, and
 I.Z. Ruzsa \cite{Ruz09} and T. Tao, V. Vu \cite{TaV06} for thorough surveys.

Unfortunately the known analogues are not
as simple in their form as the original Brunn-Minkowski
inequality. A formula due to Gardner and Gronchi says that
 if $A$ is not contained
 in any affine subspace of $\mathbb{R}^d$ then
$$
|A+B|\geq (d!)^{-\frac1d}(|A|-d)^{\frac1d}+|B|^{\frac1d}.
$$
In this paper, we discuss
a more direct version
of the Brunn-Minkowski inequality in the plane, which would
improve Freiman's inequality if both $A$ and $B$ are two-dimensional.

In the planar case ($d=2$), a recent conjecture
by Matolcsi and Ruzsa \cite{MaR} might point to the right version of the Brunn-Minkowski
inequality. Let $A$ be a finite non-collinear point set in $\R^2$. We write ${\rm tr}(A)$ to denote the common number
of triangles in any triangulation of $[A]$
using the points of $A$ as vertices. 
If  $b_A$ and $i_A$ denote
the number of points of $A$ in $\partial [A]$
and ${\rm int}[A]$, then the Euler formula yields
\begin{equation}
\label{Eulerpoints}
{\rm tr}(A)=b_A+2i_A-2.
\end{equation}
If $\Pi$ is a polygon  whose vertices are in $\Z^2$, and
 $A=\Z^2\cap \Pi$, then Pick's theorem says that
$$
{\rm tr}(A)=2\lambda(\Pi).
$$

Now the Ruzsa-Matolcsi conjecture proposes that if $A$ and
$B$ in the plane are not collinear then
\begin{equation}
\label{ruzsabrunnconj}
{\rm tr}(A+B)^{\frac12}\geq {\rm tr}(A)^{\frac12}+{\rm tr}(B)^{\frac12}.
\end{equation}
We note that equality holds if for a polygon $\Pi$ whose vertices are in $\Z^2$, and
integers $k,m\geq 1$, we have $A=\Z^2\cap k\Pi$
and $B=\Z^2\cap m\Pi$.

In this paper, we verify \eqref{ruzsabrunnconj}
in some special cases. To present our main idea we note that if $\alpha,\beta>0$, then
\begin{equation}
\label{Jensen}
(\alpha+\beta)^2\leq 2(\alpha^2+\beta^2),
\end{equation}
with equality if and only if $\alpha=\beta$.
Thus Conjecture (\ref{ruzsabrunnconj}) follows from
\begin{equation}
\label{ruzsabrunnconjbad}
{\rm tr}(A+B)\geq 2[{\rm tr}(A)+{\rm tr}(B)].
\end{equation}
This inequality does not hold in general. For example, let $\Pi$ be a polygon  whose vertices are in $\Z^2$, and
let $A=\Z^2\cap k\Pi$
and $B=\Z^2\cap m\Pi$ for 
integers $k,m\geq 1$. If $k\neq m$, then we have equality in the Brunn-Minkowski theorem for
$X=[A]$ and $Y=[B]$. Still, as we verify, (\ref{ruzsabrunnconjbad}) holds in several interesting cases.

The triangulation conjecture (\ref{ruzsabrunnconj}) can be written in the following form.

\begin{conj}[Main conjecture]
\label{main-conjecture}
If $A$ and $B$ are finite non-collinear sets $\R^2$, then
$$
\sqrt{2i_{A+B}+b_{A+B}-2}\geq \sqrt{2i_A+b_A-2}+\sqrt{2i_B+b_B-2}.
$$
\end{conj}

In turn,  (\ref{ruzsabrunnconjbad}) is equivalent with
\begin{equation}
\label{ruzsabrunnconjbadib}
2i_{A+B}+b_{A+B}\geq 4i_A+4i_B+2b_A+2b_B-6.
\end{equation}

\section{Remarks on the boundary}

In the following, we need the notion of exterior normal. A vector $u$ is an exterior normal at $x_0$ to $[A]$, where $x_0\in A$, if
$$
u\cdot x_0=\max\{u\cdot x:\, x\in A\}.
$$
It immediately follows that only points in the boundary of $[A]$ will have nonzero exterior normals. It also follows that if $a+b$ is a boundary point of $[A+B]$ for $a\in A$ and $b\in B$, then an exterior unit normal $u$ at  $a+b$ to $[A+B]$ is an exterior unit normal at $a$ to $[A]$, and at  $b$ to $[B]$. We conclude the following.

\begin{lemma}
\label{sum-boundary}
If $A$ and $B$ are finite non-collinear sets in $\R^2$, and $a\in A$ and $b\in B$, then $a+b$ lies on the boundary of $[A+B]$
with nonzero exterior unit normal vector $u$ if and only if $u$ is an exterior normal to $[A]$ at $a$ and to $[B]$ at $b$.
\end{lemma} 

For a unit vector $u$, and finite set $A$, define the collinear set of points
$$
A_u=\{x\in A:\,u\cdot x=\max_{y\in A}u\cdot y\}.
$$

\begin{lemma}
\label{sum-boundary-points}
For any  finite non-collinear sets $A$ and $B$ in $\R^2$, we have 
$$
b_{A+B}\geq b_A+b_B,
$$
with equality if and only if  $|A_u|\geq 2$ and $|B_u|\geq 2$ for a unit vector $u$ imply that
$A_u$ and $B_u$ are arithmetic progressions of the same difference.
\end{lemma}
\proof For a finite collinear set $C$, let $S(C)=|C|-1$; namely, the number of segments the points of $C$ divide the line into.
Therefore if $C$ and $D$ are contained in parallel lines, then $S(C+D)\geq S(C)+S(D)$, 
with equality if and only
if either $|C|=1$, or $|D|=1$, or $C$ and $D$ are arithmetic progressions of the same difference. Applying this observation to $C=A_u$ and $D=B_u$ for each unit vector that is an exterior normal to a side of $[A+B]$ yields the lemma.
\proofbox

\section{Sums with unique representation for each point}

In this section we consider the case where representation of points in $A+B$ is unique. We say that the representation is unique when for all $x \in A+B$, if $x=a_1+b_1$ and $x=a_2+b_2$, then $a_1=a_2$ and $b_1=b_2$.
\begin{theo}
\label{unique}
If the representation of points in $A+B$ is unique, then Conjecture~\ref{main-conjecture} holds.
\end{theo}
\proof
From the previous section, we see that whether $x=a+b\in A+B$ lies on the boundary of $[A+B]$ depends only on the exterior normals of $a \in A$ and $b \in B$. So applying any transformation to $A$ or $B$ that preserves $|A+B|$, ${\rm tr}(A)$, ${\rm tr}(B)$, and the exterior normals of $A$ and $B$ will also preserve ${\rm tr}(A+B)$. Note that scalar multiplication by $\epsilon$, where $\epsilon A=\{\epsilon a : a \in A\}$ satisfies the latter three conditions immediately. Since the representation of points in $A+B$ is unique, picking $\epsilon$ so that the representation of points in $\epsilon A+B$ is also unique will satisfy the first condition.

We pick $\epsilon$ small enough so that, for fixed $b \in B$, letting $\epsilon A+b=\{a+b:a \in \epsilon A\}$, for any $x \in \epsilon A+B$, if $x \in [\epsilon A+b]$, then $x \in \epsilon A+b$. Geometrically, this amounts to shrinking $A$ to a degree such that $\epsilon A+B$ looks like a little copy of $A$ placed at each point in $B$. It follows that the representation of points in $\epsilon A + B$ is unique, and hence ${\rm tr}(\epsilon A + B) = {\rm tr}(A + B)$.

Assume without loss of generality that ${\rm tr}(A)={\rm tr}(\epsilon A)\ge {\rm tr}(B)$. We begin to draw lines between points in $\epsilon A + B$ to form a partial triangulation, which can be extended to a triangulation of $\epsilon A +B$. For each $b \in B$, draw lines on $\epsilon A + b$ that form a triangulation of that set. Then, consider a triangulation $T$ of $B$. For each $b_1,b_2 \in B$ that are connected by a line in $T$, consider $\epsilon A+b_1$ and $\epsilon A+b_2$. Pick a point $b_1^* \in \epsilon A+b_1$ that has exterior normal $b_2 - b_1$ in $[\epsilon A+b_1]$. Pick a point $b_2^* \in \epsilon A+b_2$ that has exterior normal $b_1 - b_2$ in $[\epsilon A+b_2]$. Now, in $\epsilon A + B$, draw a line between $b_1^*$ and $b_2^*$. Geometrically, we have mimicked a triangulation of $A$ at each little copy of $A$, and a triangulation of $B$ on a large scale, treating each little copy of $A$ as a point in $B$. Letting ${\rm ptr}(\epsilon A + B)$ denote the number of polygons enclosed in this partial triangulation, it follows:
\begin{equation}
{\rm tr}(A+B)={\rm tr}(\epsilon A+B) \ge {\rm ptr}(\epsilon A+B)=|B|{\rm tr}(A)+{\rm tr}(B).
\end{equation}
Conjecture~\ref{main-conjecture} then follows from
\begin{equation}
\label{conditionforunique}
\sqrt{|B|{\rm tr}(A)+{\rm tr}(B)} \ge \sqrt{{\rm tr}(A)} + \sqrt{{\rm tr}(B)}.
\end{equation}
Since $|B| \ge 3$ and ${\rm tr}(A)\ge {\rm tr}(B)$, $(|B|-2){\rm tr}(A) \ge {\rm tr}(B)$ holds, which in turn implies (\ref{conditionforunique}).
\proofbox

\section{The case $i_A=i_B=1$}

We see that Lemma~\ref{sum-boundary-points} yields that (\ref{ruzsabrunnconjbadib}), and in turn Conjecture~\ref{main-conjecture} would follow from
\begin{equation}
\label{boundary-conjecture00}
2i_{A+B}\geq 4i_A+4i_B+b_A+b_B-6,
\end{equation}
which we have already noted does not always hold. However, in the remainder of this paper we show it holds for two special cases. The proof of the first case is simple:
\begin{theo}
\label{oneintpoint}
When $i_A=i_B=1$, Conjecture~\ref{main-conjecture} holds.
\end{theo}
\proof
From Lemma \ref{sum-boundary}, it follows that if $a \in A_{int}=\{a \in A: a \in {\rm int}([A])\}$, then $a+B\subset (A+B)_{int}$. So by (\ref{linearsum}), since $i_A$ and $i_B$ are nonempty, $i_{A+B} \ge i_{A}+|B|-1$, and similarly $i_{A+B} \ge i_{B}+|A|-1$. Thus, since $|A|=i_A + b_A$ and $|B|=i_B+b_B$, we have
\begin{equation}
2i_{A+B} \ge 2i_A+2i_B+b_A+b_B-2.
\end{equation}
In the case that $i_A=i_B=1$, (\ref{boundary-conjecture00}) follows.
\proofbox

\section{The case $|A|=b_A$ and $|B|=b_B$}

We now turn to the case $|A|=b_A$ and $|B|=b_B$, or in other words, both $A$ and $B$ lie on the boundary of their convex hulls. In this case, (\ref{boundary-conjecture00}) becomes
\begin{equation}
\label{boundary-conjecture0}
2i_{A+B}\geq b_A+b_B-6.
\end{equation}

The bad news is that (\ref{boundary-conjecture0}) does not always hold. Let
\begin{eqnarray*}
\widetilde{A}&=&\{(0,0),(1,0),(0,1)\}=\{(x,y)\in\N^2:\,x+y\leq 1\};\\
\widetilde{B}&=&\{(0,0),(1,0),(0,1),(2,0),(1,1),(0,2)\}
=\{(x,y)\in\N^2:\,x+y\leq 2\}.
\end{eqnarray*}
Therefore $|\widetilde{A}|=b_{\widetilde{A}}=3$, $|\widetilde{B}|=b_{\widetilde{B}}=6$, and
 $\widetilde{A}+\widetilde{B}=\{(x,y)\in\N^2:\,x+y\leq 3\}$ yields
$i_{\widetilde{A}+\widetilde{B}}=1$. In particular, (\ref{boundary-conjecture0}) fails to hold for $\widetilde{A}$ and 
$\widetilde{B}$, but the good news is that  Conjecture~\ref{main-conjecture} does hold for them. 

We note that  $\widetilde{B}=\widetilde{A}+\widetilde{A}$. Actually, if $A$ is any set of three non-collinear points, and $B=A+A$, then there exists a linear transform $\varphi$  such that $A$ is a translate of $\varphi\widetilde{A}$, and $B$ is a translate of $\varphi\widetilde{B}$. Therefore (\ref{boundary-conjecture0}) does not  hold for that $A$ and $B$, as well. However, in the remainder of the paper, we prove the following theorem. From this result Conjecture~\ref{main-conjecture} holds for the case when $|A|=b_A$ and $|B|=b_B$.

\begin{theo}
\label{work-conjecture}
If $A$ and $B$ are finite non-collinear sets in $\R^2$ such that $|A|=b_A$,  
$|B|=b_B$ and (\ref{boundary-conjecture0}) fails to hold,  then
either $|A|=3$, and $B$ is a translate of $A+A$, or $|B|=3$, and $A$ is a translate of $B+B$.
\end{theo}

To prove Theorem~\ref{work-conjecture}, we consider a unit vector $v$ not parallel to any side of $[A]$ or $[B]$. 
We think of $v$ as pointing vertically upwards. Let $l_{v,A}$ and $r_{v,A}$ be the leftmost and rightmost vertices of $[A]$, repectively. We note that $l_{v,A}$ and $r_{v,A}$ are unique, because $v$ is not parallel to any side of $[A]$. Similarly, let $l_{v,B}$ and $r_{v,B}$ be the (unique) leftmost and rightmost vertices of $[B]$, repectively.

Remember that $v$ points upwards. We observe that  $l_{v,A}$ and $r_{v,A}$ divide the boundary of $[A]$ into an "upper" polygonal arc, and a "lower" polygonal arc. Let $A_{v,{\rm upp}}$ and $A_{v,{\rm low}}$ denote the set of points of $A$ in the upper polygonal arc, and in the lower polygonal arc,  respectively, exluding $l_{v,A}$ and $r_{v,A}$. For $a\in A$, we have
\begin{align}
\label{Aupp}
a\in A_{v,{\rm upp}}&\mbox{ \  iff $u\cdot v>0$ for any unit exterior normal $u$ to $[A]$ at $a$};\\
\label{Alow}
a\in A_{v,{\rm low}}&\mbox{ \ iff $u\cdot v<0$ for any unit exterior normal $u$ to $[A]$ at $a$}.
\end{align}
In addition, as  $l_{v,A}$ and $r_{v,A}$ are excluded, we have
\begin{equation}
\label{Aupplow}
|A_{v,{\rm upp}}|+|A_{v,{\rm low}}|=b_A-2.
\end{equation}
Similarly, $l_{v,B}$ and $r_{v,B}$ divide the boundary of $[B]$ into an ``upper'' polygonal arc, and a ``lower'' polygonal arc, and 
$B_{v,{\rm upp}}$ and $B_{v,{\rm low}}$ denote the set of points of $B$ in the upper polygonal arc, and in the lower polygonal arc,  respectively, exluding $l_{v,B}$ and $r_{v,B}$. For $b\in B$, we have
\begin{align}
\label{Bupp}
b\in B_{v,{\rm upp}}&\mbox{ \  iff $u\cdot v>0$ for any unit exterior normal $u$ to $[B]$ at $b$};\\
\label{Blow}
b\in B_{v,{\rm low}}&\mbox{ \ iff $u\cdot v<0$ for any unit exterior normal $u$ to $[B]$ at $b$};\\
\label{Bupplow}
&|B_{v,{\rm upp}}|+|B_{v,{\rm low}}|=b_B-2.
\end{align}


\begin{lemma}
\label{upperlower}
Let $A$ and $B$ be finite non-collinear sets in $\R^2$, and let $v$ be a  unit vector  not parallel to any side of $[A]$ or $[B]$.
If  $A_{v,{\rm upp}}$, $A_{v,{\rm low}}$, $B_{v,{\rm upp}}$ and $B_{v,{\rm low}}$ are all non-empty, then
(\ref{boundary-conjecture0}) holds.
\end{lemma} 
\proof Lemma~\ref{sum-boundary}, (\ref{Aupp}) and (\ref{Blow}) yield that 
$A_{v,{\rm upp}}+B_{v,{\rm low}}\subset {\rm int} [A+B]$, therefore
$$
i_{A+B}\geq |A_{v,{\rm upp}}+B_{v,{\rm low}}|\geq |A_{v,{\rm upp}}|+|B_{v,{\rm low}}|-1.
$$
In addition, Lemma~\ref{sum-boundary}, (\ref{Alow}) and (\ref{Bupp}) yield that 
$A_{v,{\rm low}}+B_{v,{\rm upp}}\subset {\rm int} [A+B]$, therefore
$$
i_{A+B}\geq |A_{v,{\rm low}}+B_{v,{\rm upp}}|\geq |A_{v,{\rm low}}|+|B_{v,{\rm upp}}|-1.
$$
We deduce from (\ref{Aupplow}) and (\ref{Bupplow}) that
$$
2i_{A+B}\geq  |A_{v,{\rm upp}}|+|B_{v,{\rm low}}|+|A_{v,{\rm low}}|+|B_{v,{\rm upp}}|-2=b_A+b_B-6.
$$
\proofbox

In other words, Lemma~\ref{upperlower} says that if (\ref{boundary-conjecture0}) does not hold, then at least one of the sets
$A_{v,{\rm upp}}$, $A_{v,{\rm low}}$, $B_{v,{\rm upp}}$ and $B_{v,{\rm low}}$ empty. We observe that replacing $v$ by $-v$ simply exchanges $A_{v,{\rm upp}}$ and $A_{v,{\rm low}}$ on the one hand, and
$B_{v,{\rm upp}}$ and $B_{v,{\rm low}}$ on the other hand. Therefore Proposition~\ref{upperlower-fine} will refine Lemma~\ref{upperlower}.
Before that, we verify another auxiliary statement where $[p,q]$ denotes the closed line segment with end points $p,q\in\R^2$.

\begin{lemma}
\label{upperlower-segment}
Let $A$ and $B$ be finite non-collinear sets in $\R^2$, and let $v$ be a  unit vector  not parallel to any side of $[A]$ or $[B]$.
If $A_{v,{\rm low}}=\emptyset$, then $i_{A+B}\geq |B_{v,{\rm upp}}|-2$, where equality yields that
$B_{v,{\rm low}}\subset [l_{v,B},r_{v,B}]$, and the segments $[l_{v,A},r_{v,A}]$ and $ [l_{v,B},r_{v,B}]$ are parallel.
\end{lemma} 
\proof We drop the reference to $v$ in the notation. After applying a linear transformation fixing $v$, we may assume that
\begin{equation}
\label{lArA}
w\cdot v=0\mbox{ \ for $w=r_A-l_A$.}
\end{equation}
We may also assume that
\begin{equation}
\label{lArA0}
 l_A\cdot v=r_A\cdot v=0.
\end{equation}
If $r_B\cdot v>l_B\cdot v$, then we reflect both $A$ and $B$ through the line $\R v$. This keeps $v$, but interchanges the roles
of $l_A$ and $r_A$ on the one hand, and the roles of $l_B$ and $r_B$ on the other hand. Therefore we may assume that
\begin{equation}
\label{lBrB}
r_B\cdot v\leq l_B\cdot v.
\end{equation}

Understanding exterior normals helps bound interior points in $[A+B]$. As 
 $A$ has some point above $[l_A,r_A]$ by $A_{{\rm low}}=\emptyset$, 
(\ref{lArA}) yields that
\begin{align}
\label{urA}
\mbox{either $u\cdot w> 0$ or $u=-v$}&\mbox{ for any exterior unit normal $u$ at $r_A$ to $[A]$},\\
\label{ulA}
\mbox{either $u\cdot w< 0$ or $u=-v$}&\mbox{ for any exterior unit normal $u$ at $l_A$ to $[A]$}.
\end{align}
We may assume that $B_{{\rm upp}}\neq\emptyset$ (otherwise Lemma~\ref{upperlower-segment} trivially holds).
We subdivide $B_{{\rm upp}}$ into the sets
\begin{align}
\label{B-}
B_{{\rm upp}}^-=&\{b\in B_{{\rm upp}}:\,u\cdot w<0\mbox{ for any exterior unit normal $u$ at $b$ to $[B]$}\}, \\
\label{B+}
B_{{\rm upp}}^+=&\{b\in B_{{\rm upp}}:\,u\cdot w>0\mbox{ for any exterior unit normal $u$ at $b$ to $[B]$}\}, \\
\label{B0}
B_{{\rm upp}}^0=&\{b\in B_{{\rm upp}}:\,\mbox{ $v$ is an exterior unit normal $u$ at $b$ to $[B]$}\}.
\end{align}
Since for any $b\in B$, the set of all exterior unit normals $u$ at $b$ to  $[B]$ is an arc of the unit circle, the sets 
$B_{{\rm upp}}^-$, $B_{{\rm upp}}^+$ and $B_{{\rm upp}}^0$ are pairwise disjoint, and their union is $B_{{\rm upp}}$.
In addition, we define 
\begin{equation}
\label{tildeB-}
\widetilde{B}_{{\rm upp}}^-=\left\{
\begin{array}{rl}
\{l_B\}\cup B_{{\rm upp}}^- & \mbox{ if there exists $b\in B$ with  $b\cdot v< l_B\cdot v$,}\\[1ex]
B_{{\rm upp}}^- & \mbox{ if $b\cdot v\geq l_B\cdot v$ for all $b\in B$.}
\end{array}
\right.
\end{equation}
It follows that if $b\in \widetilde{B}_{{\rm upp}}^-$, then
\begin{equation}
\label{tildeB-u}
\mbox{either $u\cdot w< 0$ or $u=v$ for an exterior unit normal $u$ to $[B]$ at $b$.}
\end{equation}
Turning to $B_{{\rm upp}}^0$,  if $B_{{\rm upp}}^0\neq\emptyset$, then there exist $l_B^0,r_B^0\in B_{{\rm upp}}^0$ such
that $r_B^0-l_B^0=s w$ for $s\geq 0$, and
\begin{align}
\label{B01}
 B_{{\rm upp}}^0=&B\cap[l_B^0,r_B^0],\\
\label{B02}
v\cdot b_0=&\max\{v\cdot b:\, b\in B\}=H\mbox{ \ for all $b_0\in B_{{\rm upp}}^0$}.
\end{align}

To estimate $i_{A+B}$, we deduce  from Lemma~\ref{sum-boundary}, and from (\ref{urA}) and (\ref{tildeB-u})  on the one hand, from (\ref{ulA}) and (\ref{B+})  on the other hand, that 
\begin{equation}
\label{minusplus}
\begin{array}{rl}
r_A+\widetilde{B}_{{\rm upp}}^-\subset & {\rm int} [A+B]\mbox{ \ if }B_{{\rm upp}}^-\neq\emptyset\\[1ex]
l_A+B_{{\rm upp}}^+\subset & {\rm int} [A+B]\mbox{ \ if }B_{{\rm upp}}^+\neq\emptyset.
\end{array}
\end{equation}
We claim that if $\widetilde{B}_{{\rm upp}}^-\neq\emptyset$ and $B_{{\rm upp}}^+\neq\emptyset$, then
\begin{equation}
\label{intersection}
\left|(r_A+\widetilde{B}_{{\rm upp}}^-)\cap(l_A+B_{{\rm upp}}^+)\right|\leq 1.
\end{equation}
We observe that $r_A+x=l_A+y$ if and only if $y-x=w$, and hence $x\cdot v=y\cdot v$. However, if 
$x_1,x_2\in \widetilde{B}_{{\rm upp}}^-$
and $y_1,y_2\in B_{{\rm upp}}^+$ with $x_1\cdot v=y_1\cdot v<x_2\cdot v=y_2\cdot v$, then
$(y_2-x_2)\cdot w<(y_1-x_1)\cdot w$, which in turn yields (\ref{intersection}). We conclude by (\ref{lArA0}), (\ref{B02}),
(\ref{minusplus}) and  (\ref{intersection}) then
\begin{equation}
\label{below}
\left|\left\{z\in (A+B)\cap {\rm int} [A+B]:\, z\cdot v<H\right\}\right|\geq |\widetilde{B}_{{\rm upp}}^-|+|B_{{\rm upp}}^+|-1.
\end{equation}

We recall that there exists some $p\in A_{{\rm upp}}$, and hence $p\cdot v>0$ by $l_A\cdot v=0$.
Thus  if $B_{{\rm upp}}^0\neq\emptyset$, and $z\in \{l_A,r_A\}+B_{{\rm upp}}^0$ is different
from $l_A+l_B^0$ and $r_A+r_B^0$, then these two points of $A+B$ lie left and right from $z$.
 Since $(l_A+l_B)\cdot v<z\cdot v$, and $(p+l_B^0)\cdot v>z\cdot v$, we have
$z\in {\rm int} [A+B]$. In particular, $|\{l_A,r_A\}+B_{{\rm upp}}^0|\geq |B_{{\rm upp}}^0|+1$ yields that
\begin{equation}
\label{top}
\left|\left\{z\in (A+B)\cap {\rm int} [A+B]:\, z\cdot v=H\right\}\right|\geq |B_{{\rm upp}}^0|-1.
\end{equation}
Adding (\ref{below}) and (\ref{top}) implies $i_{A+B}\geq |B_{{\rm upp}}|-2$. If $i_{A+B}= |B_{{\rm upp}}|-2$, 
then $\widetilde{B}_{{\rm upp}}^-=B_{{\rm upp}}^-$, and hence $r_B\cdot v=l_B\cdot v$ by (\ref{lBrB}) and (\ref{tildeB-}),
and $B_{{\rm low}}\subset [l_B,r_B]$. In particular (\ref{lArA}) implies that 
$[l_{v,A},r_{v,A}]$ and $ [l_{v,B},r_{v,B}]$ are parallel.
\proofbox

\begin{prop}
\label{upperlower-fine}
Let $A$ and $B$ be finite non-collinear sets in $\R^2$, and let $v$ be a  unit vector  not parallel to any side of $[A]$ or $[B]$.
If  (\ref{boundary-conjecture0}) does not hold,  then possibly after exchanging $A$ and $B$, or $v$ by $-v$, we have the following.
\begin{description}
\item{(i)} $A_{v,{\rm low}}=\emptyset$;
\item{(ii)} $B_{v,{\rm low}}\subset [l_{v,B},r_{v,B}]$;
\item{(iii)} $[l_{v,A},r_{v,A}]$ and $[l_{v,B},r_{v,B}]$ are parallel;
\item{(iv)} either $B_{v,{\rm low}}=\emptyset$ and $b_B=b_A$, or $|B_{v,{\rm upp}}|=|A_{v,{\rm upp}}|+|B_{v,{\rm low}}|+1$
and $b_B>b_A$.
\end{description}
\end{prop} 
\proof We drop the reference to $v$ in the notation.
To present the argument, we make some preparations. Again using that  (\ref{boundary-conjecture0}) does not hold,
 Lemma~\ref{upperlower} yields that
possibly after exchanging $A$ and $B$, or $v$ by $-v$, we may assume that 
$$
A_{{\rm low}}=\emptyset.
$$
Possibly after exchanging $A$ and $B$ again, we may assume that
\begin{equation}
\label{Bbig}
\mbox{if $B_{{\rm low}}=\emptyset$, then $b_B\geq b_A$}.
\end{equation}

Since (\ref{boundary-conjecture0}) does not hold, we have
\begin{equation}
\label{no-boundary-conjecture0}
i_{A+B}<\mbox{$\frac12$}(b_A+b_B)-3.
\end{equation}

First we show that
\begin{equation}
\label{Buppcard}
\begin{array}{rl}
\mbox{either }|B_{{\rm upp}}|=\frac{b_A+b_B}2-2,&\mbox{$B_{{\rm low}}=\emptyset$ and $b_A=b_B$,}\\[1ex]
\mbox{or }|B_{{\rm upp}}|>\frac{b_A+b_B}2-2.&
\end{array}
\end{equation}
If $B_{{\rm low}}=\emptyset$, then $b_B\geq b_A$ by (\ref{Bbig}), and hence
$$
|B_{{\rm upp}}|=b_B-2\geq \mbox{$\frac{b_A+b_B}2$}-2,
$$
with equality only if $b_A=b_B$.

If $B_{{\rm low}}\neq\emptyset$, then we use that $A_{{\rm upp}}\neq\emptyset$ by $A_{{\rm low}}=\emptyset$.
Thus  Lemma~\ref{sum-boundary}, (\ref{Aupp}) and (\ref{Blow}) yield that $A_{{\rm upp}}+B_{{\rm low}}$ lies  in the interior of $[A+B]$. Combining this fact with (\ref{Bupplow}) leads to
\begin{eqnarray}
\nonumber
i_{A+B}&\geq & |A_{{\rm upp}}+B_{{\rm low}}|\geq |A_{{\rm upp}}|+|B_{{\rm low}}|-1\\
\label{iAupp1}
&=&b_A-2+b_B-2-|B_{{\rm upp}}|-1=b_A+b_B-|B_{{\rm upp}}|-5.
\end{eqnarray}
Therefore $|B_{{\rm upp}}|>\frac{b_A+b_B}2-2$ by (\ref{no-boundary-conjecture0}), proving (\ref{Buppcard}).

It follows from  (\ref{no-boundary-conjecture0}) and (\ref{Buppcard}) that $i_{A+B}<|B_{{\rm upp}}|-1$, thus
Lemma~\ref{upperlower-segment} implies that $i_{A+B}=|B_{{\rm upp}}|-2$, and in turn
 (ii) and (iii) of Proposition~\ref{upperlower-fine} hold. To prove (iv), we deduce from (\ref{no-boundary-conjecture0})  that
$$
b_A+b_B-6>2i_{A+B}=2|B_{{\rm upp}}|-4.
$$
Therefore (\ref{Buppcard}) yields that either $B_{{\rm low}}=\emptyset$ and $b_A=b_B$, or
$$
b_A+b_B-4< 2|B_{{\rm upp}}|<b_A+b_B-2.
$$
In particular, $2|B_{{\rm upp}}|=b_A+b_B-3$ in the second case, which is 
in turn equivalent with $|B_{{\rm upp}}|=|A_{{\rm upp}}|+|B_{{\rm low}}|+1$ by 
$|A_{{\rm upp}}|=b_A-2$ and (\ref{Bupplow}). In addition,
$|B_{{\rm upp}}|=|A_{{\rm upp}}|+|B_{{\rm low}}|+1$ implies that $b_B>b_A$.
\proofbox

We have now developed enough machinery to prove Theorem~\ref{work-conjecture}. We repeat it here:

\begin{theo}
\label{boundarytheorem}
If $A$ and $B$ are finite non-collinear sets in $\R^2$ such that $|A|=b_A$,  
$|B|=b_B$ and (\ref{boundary-conjecture0}) fails to hold,  then
either $|A|=3$, and $B$ is a translate of $A+A$, or $|B|=3$, and $A$ is a translate of $B+B$.
\end{theo}
\proof
We follow Proposition~\ref{upperlower-fine}, and choose $A$, $B$, and $v$ as in that result. For each $x \in A$, we have that if $x$ lies on a corner of $[A]$, there exist vectors $v_{x,l}$ and $v_{x,r}$ such that $x=l_{v_{x,l},A}$ and $x=r_{v_{x,r},A}$. Since $A_{v,low}=\emptyset$, in the first case $r_{v_{x,l},A}=r_{v,A}$ and in the second $l_{v_{x,r},A}=l_{v,A}$. Consider one such $x\in A_{v,upp}$. By Proposition~\ref{upperlower-fine}, it follows that $A$ is a subset of the triangle $T_A$ formed by $l_{v,A}$, $r_{v,A}$, and $x$. And by the same proposition, all the sides of $[B]$ must be parallel to sides in $A$, so $B$ is a subset of some triangle $T_B$=$\phi T_A$, where $\phi$ is a composition of a transposition and scalar multiplication. Then the corners of $[B]$ are $l_{v,B}$, $r_{v,B}$, and some point $y \in B$. We define the open line segments:
\begin{align}
s_1 = & (l_{v,A},r_{v,A}) \\
s_2 = & (l_{v,A},x) \\
s_3 = & (x,r_{v,A}) \\
t_1 = & (l_{v,B},r_{v,B}) \\
t_2 = & (l_{v,B},y)  \\
t_3 = & (y,r_{v,B}).
\end{align}
Let $A_i=s_i \cap A$ and $B_i=t_i \cap B$ for $i \in \{1,2,3\}$. Note that $A_1=\emptyset$, and $s_i$ is parallel to $t_i$, yet $A_i=\emptyset$ or $B_i=\emptyset$.

Assume for contradiction that $|A| > 3$. By Proposition~\ref{upperlower-fine}, $|B_{v,upp}| \ge 2$. Thus $B_i \ne \emptyset$ for one $i\in\{2,3\}$. Assume without loss of generality that $B_3 \ne \emptyset$; then by Proposition~\ref{upperlower-fine}, $A_3 = \emptyset$ and so $A_2 \ne \emptyset$. Thus, letting $p \in A_2$, since $B_1$ and $B_3$ share no nonzero exterior normals with $p$, and since $A_2$ and $r_{v,B}$ share no nonzero exterior normals, $B_1+p, A_2+r_{v,B},B_3+p\in(A+B)_{int}$. And since $T_B=\phi T_A$, these three sets are pairwise disjoint. So 
\begin{equation}
i_{A+B} \ge |B_1+p|+|A_2+r_{v,B}|+|B_3+p|=b_A+b_B-6,
\end{equation}
and thus (\ref{boundary-conjecture0}) holds, contrary to our assumption. So $|A|=3$.

By Proposition~\ref{upperlower-fine}, we have that if $B_{v,low}=\emptyset$, then $b_A=b_B=3$. So, $2i_{A+B}\ge b_A+b_B-6=0$, and again (\ref{boundary-conjecture0}) holds. Thus, we have that $|B_{v,low}|\ge 1$, and so 
\begin{equation}
|B_{v,upp}|=|B_{v,low}|+2.
\end{equation}
That is,
\begin{equation}
|B_2|+|B_3|=|B_1|+1
\end{equation}
By the same argument, we get
\begin{align}
|B_1|+|B_2|=|B_3|+1 \\
|B_1|+|B_3|=|B_2|+1
\end{align}
It follows that $|B_1|= |B_2|=|B_3|=1$ and so $b_B=6$. 

Now, $i_{A+B}>0$, and if $i_{A+B}\ge 2$ then (\ref{boundary-conjecture0}) holds, contradicting our assumption. Assuming then that $i_{A+B}=1$, we let $b_i \in B_i$ for $i \in \{1,2,3\}$. Then we see that $x+b_1=r_{v,A}+b_2=l_{v,A}+b_3$. And since $T_B=\phi T_A$, $B$ must just be a translated version of $A+A$. And, as was mentioned in the beginning of this section, Conjecture~\ref{main-conjecture} holds for $A$ and $B$.
\proofbox

\noindent K.J.  B\"or\"oczky\\ 
Alfr\'ed R\'enyi Institute of Mathematics\\
Hungarian Academy of Sciences\\	 
1053 Budapest, Re\'altanoda u. 13-15\\	 
Hungary, and

\noindent Central European University\\
1051 Budapest, N\'ador u. 9\\
 Hungary

\medbreak
\noindent B. Hoffman\\
Lewis \& Clark College \\
615 Palatine Hill Road \\
Portland, OR 97219 \\
United States

\end{document}